\documentclass[11pt,a4paper]{amsart}
\usepackage{amsrefs}
\usepackage{amssymb}

\usepackage{mathptmx} 
\usepackage[scaled=0.90]{helvet} 
\usepackage{courier} 
\normalfont
\usepackage[T1]{fontenc}

\usepackage{marvosym}
\usepackage{graphicx}
\usepackage{epstopdf}
\usepackage{lscape}
\usepackage{pifont}

\usepackage[curve,tips]{xypic}
\SelectTips{eu}{11}
\UseTips

\renewcommand{\epsilon}{\varepsilon}
\renewcommand{\setminus}{\smallsetminus}

\newtheorem{theorem}{Theorem}[section]
\newtheorem{proposition}[theorem]{Proposition}
\newtheorem{corollary}[theorem]{Corollary}
\newtheorem{lemma}[theorem]{Lemma}

\theoremstyle{definition}

\newtheorem{definition}[theorem]{Definition}

\theoremstyle{remark}

\newcommand{\cll}{{\scriptstyle\bf L}}
\newcommand{\clh}{{\scriptstyle\bf H}}

\newcommand{\lhf}{\cll\clh\mathfrak F}
\newcommand{\hf}{\clh\mathfrak F}

\newcommand{\XX}{\mathfrak X}

\newcommand{\FF}{\mathfrak F}

\newcommand{\Q}{\mathbb Q}
\newcommand{\Z}{\mathbb Z}

\newcommand{\GL}{\operatorname{GL}}

\newcommand{\fpinfty}{{\FP}_{\infty}}

\newcommand{\FP}{\operatorname{FP}}

\newcommand{\cohom}[3]{H^{{\raise1pt\hbox{$\scriptstyle#1$}}}(#2\>\!,#3)}
\newcommand{\tatecohom}[3]%
  {\widehat H^{{\raise1pt\hbox{$\scriptstyle#1$}}}(#2\>\!,#3)}

\newcommand{\Cohom}[3]%
  {H^{{\raise1pt\hbox{$\scriptstyle#1$}}}\big(#2\>\!,#3\big)}
\newcommand{\Tatecohom}[3]%
  {\widehat H^{{\raise1pt\hbox{$\scriptstyle#1$}}}\big(#2\>\!,#3\big)}

\newcommand{\homol}[3]{H_{{\lower1pt\hbox{$\scriptstyle#1$}}}(#2\>\!,#3)}
\newcommand{\homolog}[2]{H_{{\lower1pt\hbox{$\scriptstyle#1$}}}(#2)}

\DeclareMathOperator*{\aster}{\text{\LARGE{\textasteriskcentered}}}

\title{ Groups possessing extensive hierarchical decompositions }

\author{T. Januszkiewicz}
\author{P. H. Kropholler}
\author{I. J. Leary}

\address{T. Januszkiewicz, Department of Mathematics, The Ohio State University,
100 Math Tower, 231 West 18th Avenue, Columbus, OH 43210-1174, U.S.A.
and Mathematical Instutute of the Polish Academy of Sciences.
On leave from Mathematical Department of Wroc\l aw University}
\email{tjan@math.ohio-state.edu}

\address{P. H. Kropholler, Department of Mathematics, University of Glasgow, 15 University Gardens, Glasgow G12 8QW, United Kingdom}
\email{p.h.kropholler@maths.gla.ac.uk}

\address{I. J. Leary, Department of Mathematics, The Ohio State
University,
100 Math Tower, 231 West 18th Avenue, Columbus, OH 43210-1174, U.S.A.
Current address: Heilbronn Institute, Department of Mathematics,
Royal Fort Annexe, Bristol, BS8 1TW, United Kingdom}
\email{leary@math.ohio-state.edu}

\thanks{The work of Januszkiewicz was partially supported by NSF grant
  DMS-0706259; the work of Leary was partially supported by NSF grant
  DMS-0804226 and by the Heilbronn Institute.}

\date{\today} 


\begin{document}
\maketitle

\section{Introduction}

Let $\XX$ be a class of discrete groups which is closed under isomorphism.
Following \cite{k-survey}*{Definition 3.2.1} we define classes of groups $\clh_\alpha\XX$ for each ordinal $\alpha$ in the following way. 
\begin{itemize}
\item In case $\alpha=0$ we define $\clh_0\XX$ to be $\XX$.
\item In case $\alpha$ is a successor ordinal we define $\clh_\alpha\XX$ to be the class of all groups which admit a finite dimensional contractible $G$-complex with stabilizers in $\clh_{\alpha-1}\XX$. The term $G$-complex is used in the sense of \cite{tomDieck}*{Chapter II}.
\item In case $\alpha$ is a (non-zero) limit ordinal we define $\clh_\alpha\XX$ to be 
$\bigcup_{\beta<\alpha}\clh_\beta\XX$.
\end{itemize}
The class $\clh\XX$ is defined to be the union of all the classes $\clh_\alpha\XX$ as $\alpha$ runs through all ordinals. In the special case when $\XX$ is taken to be the class $\FF$ of all finite groups it has been shown by the second author and others that the class $\hf$ enjoys many interesting properties.
\begin{itemize}
\item $\hf$ is a large class of groups: see the discussion in \cite{k-survey}.
\item Not all groups belong to $\hf$: until recently the key known example of a non-$\hf$ group was Thompson's group $F$, known to be of type $\fpinfty$, \cite{browngeoghegan} and yet having infinite cohomological dimension and so failing the key finiteness theorem for $\hf$-groups of type $\fpinfty$ established in \cite{k-fpinfty}.
\end{itemize}
In many situations it is convenient to work with the class of locally $\hf$-groups. These are groups whose finitely generated subgroups belong to $\hf$ and we write $\lhf$ for this class. In general, for a class of groups $\XX$, we write $\cll\XX$ for the class of groups whose finite subsets are contained in $\XX$-subgroups.

In recent work \cite{abjlms} new examples of groups which lie outside $\hf$ have been constructed. Infinite groups $G$ are constructed which have the property that
every finite dimensional contractible $G$-complex has a global fixed point.
Such a group $G$ clearly cannot belong to $\hf$. Inspired by certain of the methods in \cite{abjlms} we shall construct new examples of groups which belong to $\hf$ and which illustrate the hierarchical nature of $\hf$. In particular we shall prove the following result, in which $\omega_1$ denotes the first uncountable ordinal.
\begin{theorem}\label{mainf}\ 
\begin{enumerate}
\item $\clh_\alpha\FF<\clh_{\alpha+1}\FF$ for every ordinal 
$\alpha\le\omega_1$.
\item $\cll\clh_{\omega_1}\FF=\lhf$.
\item $\cll\clh_\alpha\FF<\lhf$ for every ordinal $\alpha<\omega_1$.
\end{enumerate}
\end{theorem}

This thorem is a major advance on the previously known result that $\clh_2\FF<\clh\FF$. The groups constructed to prove Theorem \ref{mainf}(i) contain torsion and it remains an open problem whether there are torsion-free $\hf$-groups which are not in $\clh_3\FF$.

Before proceeding to the proof of Theorem \ref{mainf} we digress to discuss one of the known torsion-free examples of a group which has been used to show that $\clh_2\FF<\clh_3\FF$. This group was first described in \cite{k-survey2} and is there shown to belong to $\clh_4\FF\setminus\clh_2\FF$. While this implies that {\em there exists} a group in $\clh_3\FF\setminus\clh_2\FF$ it does not establish whether the constructed group is is in $\clh_3\FF\smallsetminus\clh_2\FF$ or in $\clh_4\FF\smallsetminus\clh_3\FF$. Moreover the arguments given in \cite{k-survey2} are flawed because they rely on the validity of the Bieri--Groves conjecture \cite{bierigroves} in cases which had not been established when the survey \cite{k-survey2} was written. The purpose of Section 2 here is to correct these shortcomings. Sections 3 and 4 are concerned with the proof of Theorem \ref{mainf}.

\section{Review of the earlier results that $\clh_0\FF<\clh_1\FF<\clh_2\FF<\clh_3\FF$}

The class $\clh_1\FF$ includes all groups of finite virtual cohomological dimension. The torsion-free groups in $\clh_1\FF$ are precisely the groups of finite cohomological dimension. The strict inclusion $\clh_0\FF<\clh_1\FF$ is clear from the fact that there are infinite $\clh_1\FF$-groups such as the infinite cyclic group, or more generally, any free abelian group of finite rank. Any torsion-free group of infinite cohomological dimension does not belong to $\clh_1\FF$. On the other hand, any countable directed union of groups in $\clh_\alpha\FF$ admits an action on a tree with stabilizers conjugate to groups in the directed system and therefore such groups belong to $\clh_{\alpha+1}\FF$, (see \cite{k-survey}*{Lemma 3.2.3} and Lemma \ref{tree} below). Thus the free abelian group of countably infinite rank belongs to $\clh_2\FF\setminus\clh_1\FF$ and we have $\clh_1\FF<\clh_2\FF$.

Turning to the strict inclusion $\clh_2\FF<\clh_3\FF$, only the following two examples were known prior to using the new methods of \cite{abjlms}.
\begin{itemize}
\item The free abelian group of rank $\aleph_\omega$ belongs to $\clh_3\FF\setminus\clh_2\FF$, \cite{dklt}.
\item For any transcendental real number $t$ the subgroup $G_t$ of $\GL_2(\Q(t))$ comprising matrices of the form 
\[\left(
   \begin{matrix} 
      a & b \\
      0 & 1 \\
   \end{matrix}
\right)\]
with $a>0$ belongs to $\clh_4\FF\setminus\clh_2\FF$. This is a group of orientation preserving affine transformations of the real line. 
\end{itemize}

The proof that $G_t$ does not belong to $\clh_2\FF$ rests on \cite{k-survey2}*{Lemmas 3.2.1 and 3.2.3} but depends in a non-trivial way on the fact that for each natural number $n\ge1$, there exists a subgroup $H_n\le G_t$ such that $H_n$ is of type $\FP_n$ with derived subgroup free abelian of infinite rank. A proof that $G_t$ belongs to $\clh_4\FF$ is given in \cite{k-survey2}*{Lemma 3.2.2}; here we establish that $G_t$ belongs to $\clh_3\FF$. 

The existence of metabelian groups $H_n$ is suggested by the Bieri--Groves conjecture \cite{bierigroves}, and is assumed in the first sentence of the proof of \cite{k-survey2}*{Lemma 3.2.3}. However the Bieri--Groves conjecture is known only in special cases and in fact it is only more recently that a special case sufficiently powerful to fulfil these needs has been proved. It is a consequence of the following theorem.

\begin{theorem}[Groves--Kochloukova, \cite{gk}*{Theorem 5}]
Let $Q$ be a finitely generated free abelian group and $A$ a finitely generated  (right) $\Z Q$-module. Assume that the action of $\Z Q$ on $A$ factors through an action of a quotient  $M=M_1\otimes\dots\otimes M_\ell$ of $\Z Q$, where $Q=Q_1\times\dots\times Q_\ell$ and $M_i=\Z Q_i/I_i$ is a cyclic $\Z Q_i$-module and $Q_i$ is a free abelian group with basis $(q_{i,j};\ 1\le j\le z_i)$ and $I_i$ is generated as ideal by $\{q_{i,j}-f_{i,j}|1\le j\le z_i\}$, where for every $i$ the set $\{f_{i,j}|1\le j\le z_i\}$ contains irreducible non-constant monic polynomials in $\Z[q_{i,1}]$ which are pairwise coprime in the sense that no two lie in a proper ideal of $\Z[q_{i,1}]$ and $f_{i,1}=q_{i,1}$. Assume further that $A$ is free as $M$-module. Then the split extension $G$ of $A$ by $Q$ is of type $\FP_m$ where $m=\min\{rk(Q_i)|1\le i\le\ell\}$.
\end{theorem}

The second author is indebted to Desi Kochloukova for useful conversations and for outlining the following method of applying the Groves--Kochloukova theorem above.
For a natural number $k$, let $p_k(t)$ denote the polynomial 
$\frac{t^k-1}{t-1}.$ 
\begin{lemma}\label{cyclotomicpolys}
If $k<\ell$ are coprime natural numbers then $p_k(t)$ and $p_\ell(t)$ are coprime when viewed as elements of the Laurent polynomial ring $\Z[t,t^{-1}]$.
\end{lemma}
\begin{proof}
If $k=1$ then $p_k(t)=1$ and there is nothing to prove. If $k>1$ then $\ell=km+r$ where $m,r$ are natural numbers, $k,r$ are coprime and $r<k$. Arguing inductively, there exist $f(t),g(t)\in\Z[t,t^{-1}]$ such that $$f(t)p_r(t)+g(t)p_k(t)=1.$$ We also have
$$p_\ell(t)=t^rp_m(t^k)p_k(t)+p_r(t).$$
Thus $\left[g(t)-t^rp_m(t^k)\right]p_k(t)+\left[f(t)\right]p_\ell(t)=1$.
\end{proof}

\begin{lemma}\label{desi}
Let $t$ be a fixed transcendental real number.
For each natural number $n$ there is a subgroup $H\le G_t$ of type $\FP_n$ which has derived subgroup of infinite rank.
\end{lemma}
\begin{proof}
Using Lemma \ref{cyclotomicpolys}, choose $f_1,f_2,\dots,f_{n}\in\Z[t]$ all of degree $\ge1$ and pairwise coprime in the Laurent polynomial ring $\Z[t,t^{-1}]$. We may arrange the choice so that $f_1=t$.
Let $Q$ be a free abelian group of rank $n$ with basis $q_1,\dots,q_n$. We define an action of $Q$ on $\Q(t)$ by $q_i\cdot f(t)=f_i(t)f(t)$. Let $A$ be the $\Z Q$-submodule of $\Q(t)$ generated by $1$. Then the Groves--Kochloukova Theorem shows that the split extension of $A$ by $Q$ is of type $\FP_n$. Moreover, this split extension is a subgroup of $G_t$ and the subgroup $A$ has infinite rank as an abelian group.
\end{proof}

\begin{proposition}\label{notclh2}
The group $G_t$ does not belong to $\clh_2\FF$.
\end{proposition}
\begin{proof}
This now follows from \cite{k-survey2}*{Lemmas 3.2.1 and 3.2.3} together with Lemma \ref{desi}.
\end{proof}
The remaining Lemmas in this section are concerned with proving that $G_t$ belongs to $\clh_3\FF$.
\begin{lemma}\label{tree}
Every countable group admits an action on a tree with finitely generated vertex and edge stabilizers.
\end{lemma}
\begin{proof}
Let $g_0,g_1,g_2,\dots$ be an enumeration of the elements of $G$ and set $G_i=\langle g_0,\dots,g_i\rangle$ for each $i\ge0$. Then $G$ is the union of the chain $G_0\le G_1\le G_2\le \dots$ and \cite{k-survey}*{Lemma 3.2.3} can be applied.
\end{proof}

\begin{lemma}\label{thickening}
Let $G$ be a group and $X$ a finite dimensional contractible
$G$-complex with countable stabilizers.  Then there is a finite
dimensional contractible $G$-complex $Y$ which admits a $G$-map
$Y\to X$ such that stabilizers in $Y$ are contained in 
finitely generated subgroups of $G$.  
\end{lemma}
\begin{proof}
We use an argument similar to one in~\cite{krophollermislin}.  
By the simplicial approximation theorem, there is a simplicial
$G$-complex $X'$ of the same dimension as $X$ and a $G$-equivariant 
homotopy equivalence $s:X'\rightarrow X$.  The $G$-complex $Y$ will 
be constructed as a simplicial $G$-complex, and a simplicial $G$-map 
$f:Y\rightarrow X'$ will be constructed.  

Choose a set $V$ of $G$-orbit 
representatives of $0$-simplices in $X'$.  For $v\in V$, let $G_v$
denote the stabilizer of $v$, and let $G\cdot v$ denote the orbit of
$v$.  For each $v$, let $T_v$ be a $G_v$-tree 
with finitely generated stabilizers as in Lemma~\ref{tree}, and let
$Y_v$ be the induced $G$-complex $Y_v= G\times_{G_v} T_v$.  Thus $Y_v$ 
is a 1-dimensional simplicial $G$-complex with finitely generated 
stabilizers, and there is a $G$-map $f: Y_v\rightarrow G\cdot v$ 
which is a homotopy equivalence.  Now let $Y^0$ be the disjoint union
of the subcomplexes $Y_v$, $Y^0=\bigcup_{v\in V} Y_v$, so that $f:Y^0
\rightarrow X'^0$ is a homotopy equivalence, where $X'^0$ denotes the 
0-skeleton of $X'$.  

For each $n>0$ and each $n$-simplex $\sigma=(x_0,\ldots,x_n)$ of $X'$,
define a simplicial complex $Y(\sigma)$ as the multiple join: 
$$Y(\sigma)= f^{-1}(x_0)*f^{-1}(x_1)*\cdots *f^{-1}(x_n).$$  Note that
every vertex of $Y(\sigma)$ is already contained in $Y^0$, and that
each $Y(\sigma)$ is contractible.  The map $f$ already defined extends
uniquely to a simplicial map $f:Y(\sigma)\to \sigma$, and $Y(\tau)$
is a subcomplex of $Y(\sigma)$ whenever $\tau$ is
a face of $\sigma$.  

Now define $Y$ and $f:Y\rightarrow X'$ by taking
the direct limit (indexed by the simplices of $X$) of the subspaces
$Y(\sigma)$.  For any $\sigma$, note that $f^{-1}(\sigma)= Y(\sigma)$.
The $G$-action on $Y^0$, which contains the vertex set of $Y$, extends
uniquely to a $G$-action on $Y$, and for this action $f:Y\to X'$ is
$G$-equivariant.  Since each $Y(\sigma)$ is contractible, it follows 
that $f$ is a homotopy equivalence, and hence $Y$ is contractible.  
Each vertex stabilizer in $Y$ is a finitely generated subgroup of $G$,
and so the stabilizer of any simplex of $Y$ is contained in a finitely
generated subgroup of $G$ as required.  
\end{proof} 

\begin{proposition}\label{sol}
Every finitely generated soluble group of derived length $d$ belongs to $\clh_d\FF$.
\end{proposition}
\begin{proof} For a soluble group $G$,
let $d(G)$ denote the derived length of $G$ and let $d'(G)$ denote the minimum of the derived lengths of all subgroups of finite index in $G$. We prove the formally stronger assertion that $G$ belongs to $\clh_{d'(G)}\FF$ by induction: the stated result follows because $d'(G)\le d(G)$. Note also that if $H$ is a finite extension of $K$ then $d'(H)=d'(K)$. 

If $d'(G)=0$ then $G\in\FF=\clh_0\FF$ and we are done. Suppose that $d'(G)>0$. Choose a normal subgroup $N$ of finite index in $G$ such that $d(N)=d'(G)$. The derived subgroup $[N,N]$ of $N$ is characteristic in $N$ and so normal in $G$. The quotient group $G/[N,N]$ is finitely generated abelian-by-finite and so it admits an action on a Euclidean space with finite stabilizers. Through the natural map $G\to G/[N,N]$ we therefore have an action of $G$ on a Euclidean space whose stabilizers are finite extensions of $[N,N]$. Using Lemma \ref{thickening} we can thicken this Euclidean space to a contractible finite dimensional $G$-complex in which each stabilizer $H$ is a subgroup of a finitely generated subgroup $K$ such that $[N,N]$ has finite index in $K[N,N]$. In this situation $d'(K)\le d'([N,N])\le d'(G)-1$ and it follows by induction that $K\in\clh_{d'(G)-1}\FF$ and therefore so is $H$. Hence $G$ belongs to $\clh_{d'(G)}\FF$ as required.
\end{proof}

\begin{corollary}\label{clh3}
The group $G_t$ belongs to $\clh_3\FF$.
\end{corollary}
\begin{proof}
The group $G_t$ is countable and metabelian. By Lemma \ref{tree} it acts on a tree with finitely generated stabilizers. Since these stabilizers belong to $\clh_2\FF$ by Proposition \ref{sol} with $d=2$ it follows that $G_t$ belongs to $\clh_3\FF$.
\end{proof}

\section{SQ-universality and variations on the Rips complex for hyperbolic groups}

We begin with an observation on the classical Higman--Neumann--Neumann embedding theorem.

\begin{lemma}\label{hnn}
Every countable group $\hf$-group $H$ can be embedded in a $2$-generator $\hf$-group $\bar H$.
\end{lemma}
\begin{proof}
It is well known that countable groups can be embedded in $2$-generator groups. 
In our context we need to be sure that membership of $\hf$ is preserved in the process.
The original proof in the classic paper \cite{hnn} of Higman--Neumann--Neumann does exactly this.
\end{proof}

Lemma \ref{hnn} says that every group is isomorphic to a subgroup of a quotient of a 
free group of rank $2$, in other words that the free group on $2$ generators is SQ-universal,
and moreover the embedding has desirable properties with respect to the $\hf$ class.
Recently it has been shown by Olshanskii \cite{o} that non-elementary hyperbolic 
groups are SQ-universal.

For our purposes we need to control $\hf$ membership just as in the case of Lemma \ref{hnn}.
This comes from combining results of Arzhantseva--Minasyan--Osin \cite{amo} with 
those of Dahmani \cite{d}.
The results of these authors that we need involve the notion of relative hyperbolicity for groups. This idea, formalised by Bowditch in \cite{bowditch}*{Section 4} and by Farb in \cite{farb}*{Definition 3.1}, is intended to enable certain groups which act on hyperbolic spaces but which are not hyperbolic in Gromov's sense to be treated within the theory of Gromov hyperbolic groups. An example to consider is the fundamental group of a hyperbolic knot complement in $S^3$ which is hyperbolic relative to the $\Z^2$ subgroup determined by a choice of basepoint on the boundary. Associated to a hyperbolic group and a specified finite generating set there are Rips complexes determined by a single parameter $d$ and the fundamental theorem \cite{salem}*{Th\'eor\`eme 12} of Rips states that for sufficiently large $d$ these are contractible: they are always finite dimensional provided $d<\infty$. Meintrup and Schick give a proof \cite{meintrupschick} that for yet larger values of the parameter these Rips complexes are classifying spaces for proper actions. More precisely Rips theorem is known to hold if $d\ge4\delta+2$ where $\delta$ is the hyperbolicity constant, and the Meintrup--Schick argument works for $d\ge16\delta+8$. The theorem of Dahmani concerns the existence of relative Rips complexes associated to relatively hyperbolic groups. The version suited to this paper is closest to \cite{d}*{Theorem 6.2} which we restate here in the following way.

\begin{theorem}\label{dahmani's 6.2}
Let $\Gamma$ be a relatively hyperbolic group in the sense of Farb,  relative to a subgroup $C$, and satisfying the bounded coset penetration property. Then, $\Gamma$ acts on a simplicial complex which is aspherical, finite dimensional, locally finite everywhere except at the vertices, with vertex stabilizers being the conjugates of $C$.
\end{theorem}

It has become standard to always include the property of bounded coset penetration within Farb's definition. The resulting definition is equivalent to a number of others, including Bowditch's \cite{bowditch}*{Section 4}. Following Bowditch, we formulate the definitions in the following way, beginning with a graph theoretic definition.

\begin{definition}\label{fine}
A {\em fine hyperbolic graph} is a graph whose geometric realisation is a geodesic metric space satisfying the condition for Gromov hyperbolicity and in addition having only finitely many circuits of length $n$ through any given edge, for all $n$. Here a {\em circuit} is a cycle which has no self-intersection.
\end{definition}

\begin{definition}\label{rhg}\cite{bowditch}*{Section 0, Definition 2}
A group $G$ is {\em relatively hyperbolic relative to a subgroup} $H$ if it admits an action on a connected fine hyperbolic graph with finite edge stabilizers, finitely many orbits of edges, and each infinite vertex stabilizer conjugate to $H$. If such a graph exists it can be chosen to have no cut-vertices, \cite{bowditch}*{Lemma 4.7}.
\end{definition}

Provided $H$ is finitely generated, the definition used by Dahmani underlying Theorem \ref{dahmani's 6.2} is equivalent to Definition \ref{rhg}. For the proof we refer the reader to \cite{szcz,bumagin,d}. For a survey of this and other equivalences between definitions see \cite{osin}*{Section 7, Appendix}.
From the point of view of Definition \ref{rhg} the relative Rips complex is elegantly described using the approach \cite{mineyevyaman} of Mineyev and Yaman. This is an alternative to the original definition which is implicit in Dahmani's \cite{d}*{Theorem 6.2} and in \cite{dahmaniyaman}*{Theorem 2.11}.

\begin{definition}\cite{mineyevyaman}*{Definition 14}
Let $d$ and $r$ be positive integer parameters.
Let $G$ be relatively hyperbolic relative to a subgroup $H$ in the sense of Definition \ref{rhg} with a graph having no cut-vertices. Then the Mineyev--Yaman relative Rips complex is the flag complex with vertices the vertices of the graph and edges the paths of length $\le d$ such the angle at any vertex of the path is $\le r$. Here the {\em angle} at a vertex of the path is defined to be the minimum length of a circuit which passes through the same two edges of the path at that vertex. Thus the Mineyev--Yaman complex is determined by the group and graph together with two parameters.
\end{definition}

For any finite values of the parameters $d$ and $r$ this complex is finite dimensional \cite{mineyevyaman}*{Corollary 17} and according to \cite{mineyevyaman}*{Theorem 19}, if $d=r$ is sufficiently large then the complex is contractible. Moreover the stabilizers of simplices of dimension $\ge1$ are finite.

\begin{proposition} \label{quot}
Let $\XX$ be a subgroup closed class of groups. 
Suppose that $H$ is a countable group in  $\XX$ 
and that $K$ is a non-elementary hyperbolic group in $\XX$. 
Then there is a quotient $Q$ of $K$ which belongs to 
$\clh\XX$ and which contains a subgroup isomorphic to $H$.
\end{proposition}

\begin{proof}
Given groups $H$, $K$, let $\bar H$ be a 2-generator group with $\bar H\ge H$
constructed using Lemma \ref{hnn}.
Results of \cite{amo} give a quotient $Q$ of $K$, which contains $\bar H$,
and which is hyperbolic relative to $\bar H$.
Now we may use a relative Rips complex as in either Dahmani or Mineyev--Yaman.  Applying this construction to $Q$
gives a finite dimensional contractible simplicial $Q$-complex in
which all infinite stabilizers are conjugate to $\bar H$.
\end{proof}

\begin{proposition}\label{build}
Let $H$ be a countable group. 
For each natural number $n$ there is a group $Q_n$ which satisfies the following conditions:
\begin{itemize}
\item $Q_n$ has a subgroup isomorphic to $H$, and
\item every contractible $Q_n$-complex of dimension $\le n$ has a global fixed point.
\end{itemize}
Moreover, if $\XX$ is a class of groups which contains all finite groups
and $H$ belongs to $\clh\XX$ then $G$ can be chosen from $\clh\XX$ also.
\end{proposition}

\begin{proof}
In Section 3 of \cite{abjlms} a group $G_{n,p}$ is constructed for each natural number 
$n$ and each prime $p$ with the following properties.
\begin{itemize}
\item $G_{n,p}$ is a non-elementary hyperbolic group;
\item every mod-$p$ acyclic $G_{n,p}$-complex of dimension $\leq n$ has a global fixed point.
\end{itemize}

For any prime $p$ the group $Q_n$ obtained from $H$ and $G_{n,p}$ using Proposition \ref{quot} has the desired properties.
\end{proof}

\section{Proof of Theorem \ref{mainf}}

\begin{theorem}\label{general}
Let $\XX$ be a subgroup closed class of groups which contains the class of all finite groups. Assume that there is a countable group which belongs to $\clh_1\XX$ but does not belong to $\XX$. Then $\clh_\alpha\XX<\clh\XX$ for all countable ordinals $\alpha$.
\end{theorem}
\begin{proof}
We shall prove by induction on $\alpha$ that
\begin{itemize}
\item there is a countable group $H_\alpha$ in $\clh\XX\setminus\clh_\alpha\XX$ for all countable ordinals $\alpha$.

\end{itemize}

If $\alpha=0$ this follows from the assumption that $\XX<\clh_1\XX$. 

If $\alpha$ is a countable limit ordinal and we have choices $H_\beta\in\clh\XX\setminus\clh_\beta\XX$ for each $\beta<\alpha$ then we can define $H_\alpha$ to be the free product
$$\aster_{\beta<\alpha}H_\beta.$$ This is a free product of countably many countable groups so is also a countable group. 
Being a free product there is an action on a tree with vertex stabilizers conjugate to the $H_\beta$ and with trivial edge stabilizers, so it follows that $H_\alpha$ belongs to $\clh\XX$. 
Since $\XX$ is subgroup closed so also is $\clh_\beta\XX$ for any $\beta$. Since $H_\alpha$ has a subgroup $H_\beta$ which does not belong to
$\clh_\beta\XX$ for each $\beta<\alpha$ it follows that
$H_\alpha$ does not belong to any of the $\clh_\beta\XX$ and therefore $H_\alpha$ is not in $\clh_\alpha\XX=\bigcup_{\beta<\alpha}\clh_\beta\XX$.

Finally if $\alpha$ is a successor ordinal then we may suppose given a group $H\in\clh\XX\setminus\clh_{\alpha-1}\XX$ by induction. For each natural number $n$ let $Q_n$ be a group satisfying all the conclusions of Proposition \ref{build}.
We can then take $H_\alpha$ to be the free product
$$\aster_nQ_n$$
as $n$ runs over the natural numbers. Since $H_\alpha$ is a free product of groups which belong to $\clh\XX$ we have that $H_\alpha$ belongs to $\clh\XX$ and it remains only to prove that $H_\alpha\notin\clh_\alpha\XX$.
Suppose by way of contradiction that $H_\alpha$ belongs to $\clh_\alpha\XX$. Then there is a finite dimensional contractible $H_\alpha$-complex $X$ with stabilizers in $\clh_{\alpha-1}\XX$. For $n>\dim X$, the subgroup $Q_n$ fixes a point of this complex and so belongs to $\clh_{\alpha-1}\XX$, but this contradicts the fact that $Q_n$ contains an isomorphic copy of the group $H$ which does not belong to $\clh_{\alpha-1}\XX$.
\end{proof}

\begin{proof}[Proof of Theorem \ref{mainf}]
By Theorem \ref{general} we have that $\clh_\alpha\FF<\clh\FF$ and therefore
$\clh_\alpha\FF<\clh_{\alpha+1}\FF$ for all countable ordinals $\alpha$. This establishes part (i) in case $\alpha<\omega_1$. 
For each countable ordinal $\alpha$, let $H_\alpha$ be a group in $\clh_{\alpha+1}\FF\setminus\clh_\alpha\FF$. Then the free product $$\aster_{\alpha<\omega_1}H_\alpha$$ does not belong to $\clh_{\omega_1}\FF$ and also, being a free product of $\clh_{\omega_1}\FF$-groups it admits an action on a tree witnessing membership of $\clh_{\omega_1+1}\FF$. This establishes (i) in the remaining case $\alpha=\omega_1$.
For part (ii) observe that every $\hf$-group of cardinality $<\kappa$ must appear in $\clh_\alpha\FF$ for some ordinal $\alpha$ of cardinality $<\kappa$. Therefore $\clh_{\omega_1}\FF$ contains all countable $\hf$-groups and in particular it contains all finitely generated $\hf$-groups. This shows that $\cll\clh_{\omega_1}\FF$ contains all $\hf$-groups and (ii) follows. Part (iii) follows from the consequence of Lemma \ref{tree} that all countable $\cll\clh_\alpha\FF$ groups belong to $\clh_{\alpha+1}\FF$.
\end{proof}


\nocite{dahmaniyaman}

\bibliography{jkl_140709}
\bibliographystyle{abbrv}

\end{document}